\documentclass[amsart]{article}
\usepackage[latin1]{inputenc}
\usepackage{amsfonts}
\usepackage{amsmath}
\usepackage{amssymb}
\usepackage[T1]{fontenc}
\usepackage{graphicx}
\usepackage{geometry}
\usepackage{mathrsfs}
\usepackage{setspace}
\usepackage{gensymb}
\usepackage{placeins}
\usepackage{color}
\usepackage{verbatim}
\usepackage[english]{babel}
\usepackage{authblk}
\usepackage{url}

\newcommand{\R}{\mathbb R}

\newcommand{\dhk}{D_h^k}

\newtheorem{remark}{Remark}
\newtheorem{definition}{Definition}
\newtheorem{proposition}{Proposition}

\begin{document}
\twocolumn
\title{Recovering Fine Details from Under-Resolved Electron Tomography Data using Higher Order Total Variation $\ell_1$ Regularization\thanks{This work is supported in part by the grants NSF-DMS 1502640 and AFOSR FA9550-15-1-0152.}}
\author[1]{Toby Sanders\thanks{toby.sanders@asu.edu}}
\author[2]{Anne Gelb}
\author[1]{Rodrigo Platte}
\author[3]{Ilke Arslan}
\author[4]{Kai Landskron}
\affil[1]{School of Mathematical and Statistical Sciences, Arizona State University}
\affil[2]{Department of Mathematics, Dartmouth College}
\affil[3]{Fundamental and Computational Sciences Directorate, Pacific Northwest National Laboratory}
\affil[4]{Department of Chemistry, Lehigh University}

\date{}
\maketitle

\begin{abstract}
 Over the last decade or so, reconstruction methods using  $\ell_1$ regularization, often categorized as compressed sensing (CS) algorithms, have significantly improved the capabilities of high fidelity imaging in electron tomography.  The most popular $\ell_1$ regularization approach within electron tomography has been total variation (TV) regularization.  In addition to reducing unwanted noise, TV regularization encourages a piecewise constant solution with sparse boundary regions.  In this paper we propose an alternative $\ell_1$ regularization approach for electron tomography based on higher order total variation (HOTV).  Like TV, the HOTV approach promotes solutions with sparse boundary regions. In smooth regions however, the solution is not limited to piecewise constant behavior.  We demonstrate that this allows for more accurate reconstruction of a broader class of images -- even those for which TV was designed for -- particularly when dealing with pragmatic tomographic sampling patterns and very fine image features.  We develop results for an electron tomography data set as well as a phantom example, and we also make comparisons with discrete tomography approaches.  
\end{abstract}

\section{Introduction}

Tomography refers to the process of non-invasive imaging via an inversion of a sequence of projections or integrals of the image.  It has long been a popular technique for 2-D and 3-D imaging across a variety of applications at many different scales.  Using the electron microscope for data acquisition, electron tomography is implemented for 3-D nanoscale image reconstruction for biological and material characterization \cite{Frank2006}.  Image formation in electron tomography requires solving an ill-posed inverse problem from limited data that can suffer from undesirable artifacts.  Thus as the computing resources and data acquisition methods evolve, research advances continue for inverse methods related to the imaging techniques for electron tomography, \cite{GorisTV,GorisTVDART,Leary,DIPS}.

Due to their simplicity, direct inversion methods such as gradient decent least squares methods and filtered backprojection are traditionally employed in tomographic applications, especially when fewer computational resources are available.  The ill-conditioning may be ameliorated by employing quadratic or Tikhonov regularization. More recently, showing far more promise for electron tomography (and many other applications) are sparsity based regularization methods that rely on minimizing an $\ell_1$ norm.  Such methods are often called $\ell_1$ regularization techniques or compressed sensing \cite{eldar2012compressed}, and they rely on the implementation of additional prior knowledge about the probable smoothness characterizations of the images.  Within the realm of $\ell_1$ techniques, electron tomography has primarly been limited to the use of the popular total variation (TV) minimization \cite{Zhang}, and it has proven overwhelmingly superior to direct inversion techniques  \cite{Leary,GorisTV}.  We note that other regularization techniques which integrate the reconstruction and segmentation of the images have also become popular in electron tomography \cite{GorisTVDART,DART,DIPS}.  These methods are collectively referred to as discrete tomography.

In this paper we present an alternative higher order total variation (HOTV) regularization approach for electron tomography.  In our investigation we use a particular form of HOTV regularization which is also called polynomial annihilation (PA) regularization \cite{Archibald2015,Archibald:2005:PFE:1061182.1068411}.  As with TV regularization, PA regularization encourages solutions with sparse boundary regions.  However, TV regularization is designed under the assumption that away from these boundaries the solution is essentially a piecewise constant function, i.e.~a polynomial of degree zero.  Alternatively, using the PA regularization assumes that the underlying solution is a piecewise polynomial of greater degree, which allows a more accurate reconstruction for a wider class of images.  Moreover, as will be demonstrated in our numerical results, using the PA regularization allows for more accurate recovery of fine features from under-resolved data.  Following up on the general framework in \cite{Archibald2015}, we look further into the performance of PA regularization on image reconstruction from electron tomography data.  We show that for piecewise constant images with many neighboring jumps,  using the PA regularization achieves more accurate approximations than using TV.  This is particularly true with typical tomographic sampling patterns in electron tomography, where the projection data are traditionally collected with parallel beam geometry at neighboring angles.

From a mathematical perspective, this work also compares with the work of \cite{VOTV}, in the observed function approximation near edges when using PA regularization.  We also note that different motivations have been used in designing similar HOTV approaches to PA regularization for imaging problems, \cite{TGV,HOTV,1257394}.

The rest of this paper is organized as follows.  In Section \ref{sec:notation} we introduce the tomographic reconstruction problem and some practical issues, including the limited sampling concerns.  In Section \ref{sec:inversemethods} we define the regularization methods and look into the appropriate parameter selection.  Numerical results are presented
in Sections 4-6, where we first consider an experi-
mental data set and then compare these results with
discrete tomography, culminating in the presentation
of some tomographic simulations. Section \ref{sec:conclusion} provides some concluding remarks.

\section{Notation and Problem Description}
\label{sec:notation}
We seek to recover images or functions denoted by $f \in \R^{n\times n}$, where the $(i,j)$ entry of $f$, denoted $f_{i,j}$, represents the image intensity at pixel $(i,j)$.  For convenience $f$ will sometimes be re-indexed as a vector, $f = \{f_i \}_{i=1}^{n^2}$.

For $p\ge 1$, the $\ell_p$ norm of $f$ is defined by 
\begin{equation}\label{lpnorm}
\| f\|_{p} = \left( \sum_{i=1}^{n^2} |f_i|^p \right)^{1/p}.
\end{equation}
For a matrix $A\in \R^{m\times n}$, the $\ell_p$ norm is defined by
\begin{equation}\label{matrix-norm}
 \| A \|_p = \max_{x\in \R^n} \frac{\| Ax \|_p}{\|x\|_p}.
\end{equation}

In electron tomography, the overarching goal is to accurately reconstruct 3-D structures from 2-D Radon transforms or projections.  We formulate the general problem for 2-D functions, and the problem is easily extended to 3-D by considering the 3-D problem as a sequence of 2-D slices\footnote{Although the 3-D problem can be considered as a sequence of 2-D problems, in practice it is beneficial to set the problem up in 3-D so that regularization of the solution can be utilized in all dimensions.}.  Thus the acquired data are values of the Radon transform of a continuous valued function $f$ are defined by
\begin{equation}\label{Radon}
R f(x, \theta) = \int_\R f((x ~ y) Q_\theta ) ~ dy ,
\end{equation}
where
\begin{equation}\label{rotation-matrix}
Q_{\theta} = 
\left[
\begin{array}{cc}
\cos \theta & \sin \theta \\
- \sin \theta & \cos \theta
\end{array}
\right].
\end{equation}

In biological and materials science, the set of values of the Radon transform given at a fixed angle $\theta$ is often called a projection of $f$, and for electron tomography {this data set} is most accurately acquired with the electron microscope in HAADF-STEM mode.  In practical terms, a projection of $f$ at angle $\theta$ can be defined as the set of all line integrals of $f$ at angle $\theta$ with respect to some coordinate system.  A projection is depicted in Figure 1.
\begin{figure}[h]
\centering
\includegraphics[scale=.15]{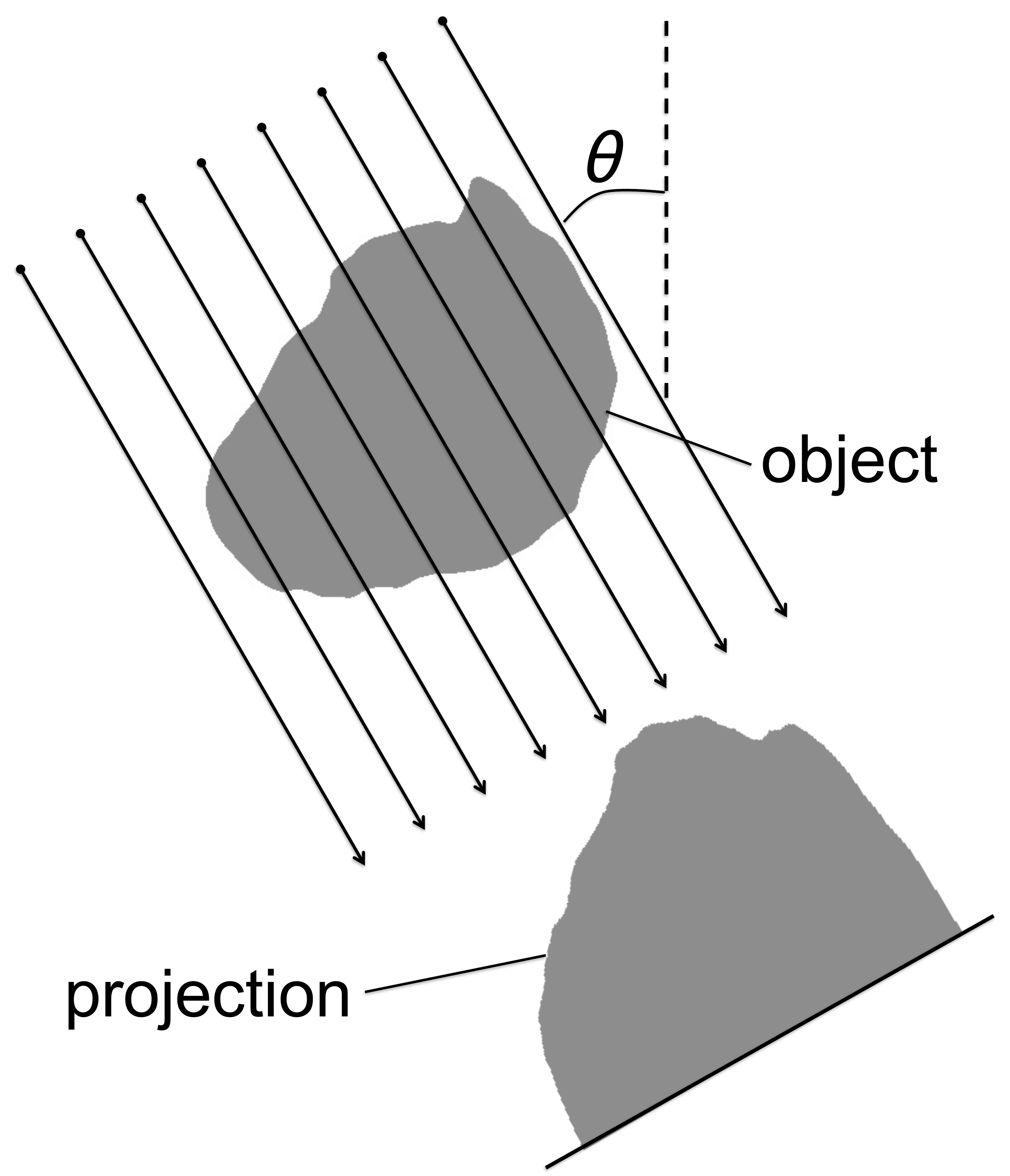}
\caption{The Radon transform of projection of a simple object at an angle of $\theta = 30\degree$.}
\end{figure}

In practice we have a finite discrete set of projection data $\{Rf(x_t , \theta_\ell) \}_{t ,\ell = 1}^{N , \alpha}$ that is corrupted with noise which we denote by
\begin{equation}
\label{eq:noisyRf}
\tilde R f_{t,\ell} = Rf(x_t , \theta_\ell) + \epsilon_{t,\ell}.
\end{equation} 
The discretization, or detection count $N$, over the $x$-coordinate is typically sufficiently small relative to the scene size, although sparse sampling strategies have been recently considered for limiting the beam dosage \cite{Saghi2015}.  The number of angles, $\alpha$, generally ranges from 50 to 100, and reducing this number has also been studied as means of reducing data acquisition time and beam dosage \cite{DIPS,DART}.  In addition, the available angular range is often limited to significantly less than the ideal $180\degree$, with a typical maximum total tilt range of $140\degree$.  This problem is infamously known as the ``missing wedge'' \cite{Arslan2006994}.

The projection values in (\ref{eq:noisyRf}) are then used to reconstruct an approximation of the original object using the linear system
\begin{equation}\label{linear-system}
Wf = b \, , \quad \text{where} \quad 
\quad 
b = 
\left(
\begin{array}{c}
\tilde Rf_{1,1}\\
\tilde Rf_{2,1} \\
 \vdots \\
\tilde Rf_{N-1,\alpha} \\
\tilde Rf_{N, \alpha}
\end{array}
\right),
\end{equation}
and $W$ is a linear operator that maps $f \in \R^{n\times n}$ to the discretized approximation of (\ref{Radon}).  {{More information on the construction of $W$ can be found in \cite{sanders2015image} on pages 8-9 within Section 1.5.}}

The direct inverse problem (attempting to solve (\ref{linear-system}) directly) \cite{Natterer} is ill-posed, and therefore will  benefit from appropriate regularization, which is discussed in Section \ref{sec:inversemethods}.  In this setting, the number of equations is usually significantly less than the desired resolution of the image, i.e. $m = \alpha N \ll n^2$.  This leads to an underdetermined system where the theory of compressive sensing may be applicable under certain assumptions, \cite{CSincoherence,DeVore2007918}.  However, even where the theory fails to hold, formulating an optimization problem that promotes sparsity of some underlying features of the image may still be effective. 

\section{Regularized Inverse Methods}
\label{sec:inversemethods}
A common model {used in determining} regularized solutions to (\ref{linear-system}) is a convex optimization problem that takes the form
\begin{equation}\label{gen-inverse}
 J (f) = \frac{\lambda}{2} \| Wf - b \|_2^2 +  H({f}),
\end{equation}
where $H( f )$ is typically some norm or semi-norm and acts as a penalty or regularization term that discourages unfavorable solutions to the problem.

Our interest in this paper is the aforementioned $\ell_1$ regularization that encourages sparsity of $f$ in some appropriate domain by setting the regularization term to $\| {Tf} \|_1 $, where $ T$ is the linear transformation under which the solution is assumed to be sparse.  Most commonly in electron tomography is to choose the regularization term as the TV norm \cite{Leary,Saghi2015,GorisTV,Sanders2015}, which is equivalent to the first order finite difference operator that maps $f$ to the differences between all adjacent pixels.

\subsection{Polynomial Annihilation Transform}
{More recently some models have used a high order TV (HOTV) approach for $\ell_1$ regularization.}  Here the linear transform $T$ can been defined as a higher order finite difference operator.  As stated previously, this operator is also known as the polynomial annihilation (PA) transform,  \cite{Archibald2015,Archibald:2005:PFE:1061182.1068411,VOTV}.  Loosely speaking, $T$ annihilates $f$ at grid points where locally $f$ is essentially a polynomial less than a desired degree, or the order of the PA operator.  We denote the PA operator of order $k$ by $T_k$ noting that $T_1$ is equivalent to TV.  As with TV, the PA regularization encourages sparse boundary regions.  However, in contrast with TV, the smooth regions are modeled as polynomials rather than constant functions.  We briefly outline the general PA methodology.  A more formal discussion of polynomial annihilation can be found in \cite{Archibald:2005:PFE:1061182.1068411}.

As outlined below, in the case of an equally spaced grid, the PA operator is  equivalent to a kth order finite difference (if the grid is not equally spaced, then the operator {takes on the more general formulation for divided differences}). Hence we will show below that the PA transform has the following definition:

\begin{definition}
Let $ f \in \R^N$ and let $j\le N - k$.  Then the $j^{th}$ element of the order $k$ PA transform of $f$ is given by
\begin{equation}\label{finite-diff}
 ({T}_k {f})_j = \sum_{m=0}^k (-1)^{k+m}  {k \choose m} f_{j+m} .
\end{equation}
\end{definition}

For indices $j>N-k$, we may also define $(T_k f)_j$ using a periodic extension of $f$.  We note that our definition differs slightly from the one presented in \cite{Archibald:2005:PFE:1061182.1068411}. It does not effect the methodology we present, however.  Observe that in multiple dimensions Definition 3.1 is simply extended along each grid direction, which we remark further on later.

To see that ${T_k}$ defined in (\ref{finite-diff}) indeed annihilates polynomials of degree less than $k$, consider first a function $f$ in the continuous domain.  Here the PA operator is evidently defined to be the $k^{th}$ derivative of $f$, since polynomials of degree less than $k$ vanish after applying $k$ derivatives.  However, if a function is defined at equally spaced grid points, then this operator simply becomes high order finite difference operator,  which is defined as (see for instance, \cite{FDBook}):
\begin{definition}[Finite-differences]\label{def2}
Let $f$ be a well defined and bounded function over $\R$.  The order $k$ finite difference operator of step size $h$, denoted $D_h^k$, applied to $f$ at $x$, is defined by
\begin{itemize}
 \item $D_h^1 f(x) = D_h f(x) = f(x+h)-f(x)$.
 \item $D_h^{k+1} f(x) = D_h^1 \left(D_h^k f(x) \right)$,  for  $k\ge1$.
\end{itemize}
\end{definition}

By definition of the derivative, for sufficiently smooth functions we have 
\begin{equation}
\lim_{h\rightarrow 0} \frac{\dhk f(x)}{h^k} = \frac{d^k}{dx^k}f(x).
\end{equation}
Additionally, repeated application the mean value theorem yields
\begin{equation}\label{mean-value}
{\dhk f(x)} = {h^k} \frac{d^k f}{dx^k}(\xi),
\end{equation}
for some $\xi \in [x,x+kh]$.  Therefore, if $f$ is a polynomial of degree less than $k$ in the interval $[x,x+kh]$, by (\ref{mean-value}) we see that $\dhk f(x) =0$.  In addition, if $f$ can be approximated reasonably well by a polynomial of degree $k$ in the interval $[x,x+kh]$, then $\dhk f(x)$ will still be very small.

Finally, we have the following closed formula for the finite differences.
\begin{definition}
\label{prop:fdformula}
 The order $k$ finite difference operator is given by the formula
 \begin{equation}
  \dhk f(x) = \sum_{m=0}^k (-1)^{m+k} {k \choose m} f(x+mh).
  \end{equation}
\end{definition}

Definition \ref{prop:fdformula} can be derived from Definition \ref{def2} by a straight forward induction over $k$, yielding exactly (\ref{finite-diff}) for discrete functions defined at equally spaced grid points.  To this end, the regularization method  can be expressed as finding the solution to
\begin{equation}\label{gen-PA}
\min_f  \Big\{ \frac{\lambda}{2}\| Wf - b \|_2^2 +  \| {T_k f} \|_1 \Big\} ,
\end{equation}
where $T_k$ is defined by (\ref{finite-diff}).

\begin{remark}
\label{rem4}
Observe that $T_1$ is equivalent to the TV regularization.  Similarly, the case $k=0$ simplifies to ${T_0 = I}$, the identity.  Finally, it is worth noting that for $k\ge1$, $T_k = T_1^k$.  Therefore, for one to move from TV to higher order PA methods, the PA transform may be defined simply by repeated application of the transform for TV.
\end{remark}

\begin{remark}
\label{rem5}
For multiple dimensional higher order finite differences, we simply extend $T_k$ to evaluate the difference along each dimension.  For instance, in 2-D, we may have 
$T_k =\left[ \begin{array}{c} T_k^x \\ T_k^y\end{array}\right]$
, where $T_k^x$ computes the differences along the $x$ dimension and likewise for $T_k^y$.  For example, for $f\in R^{n\times n}$, $(T_k^x f)_{j,\ell} = \sum_{m=0}^k (-1)^{k+m} {k\choose m} f_{j+m, \ell}$
\end{remark}

\begin{remark}
\label{rem6}
 In order to solve (\ref{gen-PA}) for large imaging problems, an efficient iterative algorithm is needed. For these general types of large $\ell_1$ optimization problems, most commonly it has been proposed to split the objective functional into equivalent subproblems, each of which can be solved relatively fast \cite{bregman,Zhang}.  In our work, we have followed the approach in \cite{Li2013}, which rewrites the problem with an augmented Lagrangian function and uses alternating direction minimization to solve the subproblems in a method sometimes called alternating direction method of multipliers (ADMM).
\end{remark}

\subsection{Regularization parameter selection}\label{lambda-select}
Generally speaking, selection of $\lambda$ depends on the signal to noise ratio (SNR).  In particular, for higher SNR one should select smaller $\lambda$ and visa versa.  However, $\lambda$ typically also depends on the relative energy of the signal, as well as properties of the sensing matrix $W$, which in this case is determined primarily by the projection geometry.  While parameter selection will be addressed more rigorously in future work, below we briefly describe how this can be done.

We observe that when solving a problem such as (\ref{gen-PA}), the appropriate selection of $\lambda$ can often be greatly simplified by first rescaling the operators $W$ and $T_k$, as this reduces the dependence on the specific problem, that is, the particular projection geometry and chosen dimension of the recontruction.  Hence $\lambda$ becomes primarily dependent on the noise levels and error in the forward model.  

The sensing matrix $W$ is rescaled so that $\|W\|_2=1$, and the data values are rescaled accordingly.  Determining this rescaling factor can be accomplished relatively easily using eigenvalue power methods (see, for instance \cite{Golub}), making this rescaling a practical rule.

Similarly, we may rescale $T_k$ so that $\|T_k\|_1=1$.  However, for simplicity we keep the integer valued coefficients defined in the operator given by (\ref{finite-diff}), and thus instead choose to rescale $\lambda$ depending on $\|T_k\|_1$.  To determine this rescaling, we offer the following proposition. 
\begin{proposition}
 For the operator $T_k$ defined in (\ref{finite-diff}), the $\ell_1$ norm of $T_k$ is given by
 \begin{equation}\label{Tk-norm1}
  \|T_k \|_1 = 2^k.
 \end{equation}
\end{proposition}
The details of this proof are given in \cite{Sanders-lambda}.
From (\ref{matrix-norm}) we have
\begin{equation}
\label{eq:T_k2k}
||T_kf||_1 \le 2^k ||f||_1,
\end{equation}
suggesting that the regularization parameter in (\ref{gen-PA}) using $T_k$ be of the form
\begin{equation}
\label{eq:lambdaselection}
\lambda = 2^{k-1} \lambda_1,
\end{equation}
for some initial choice $\lambda_1$ that is ideal for TV.  It is important to note that equality may only be achieved in (\ref{eq:T_k2k}) for signals that are not typically of interest, and indeed most signals will produce small values for $\|T_kf \|_1$.  However, the numerical results in this paper strongly support choosing $\lambda$ using (\ref{eq:lambdaselection}).  We refer the reader to \cite{Sanders-lambda} for additional theoretical results.


\section{Electron Tomography on a Mesoporous Boron Nitride Nanoparticle}\label{ET-1}

\begin{figure*}
 \centering
 \includegraphics[width=1\textwidth]{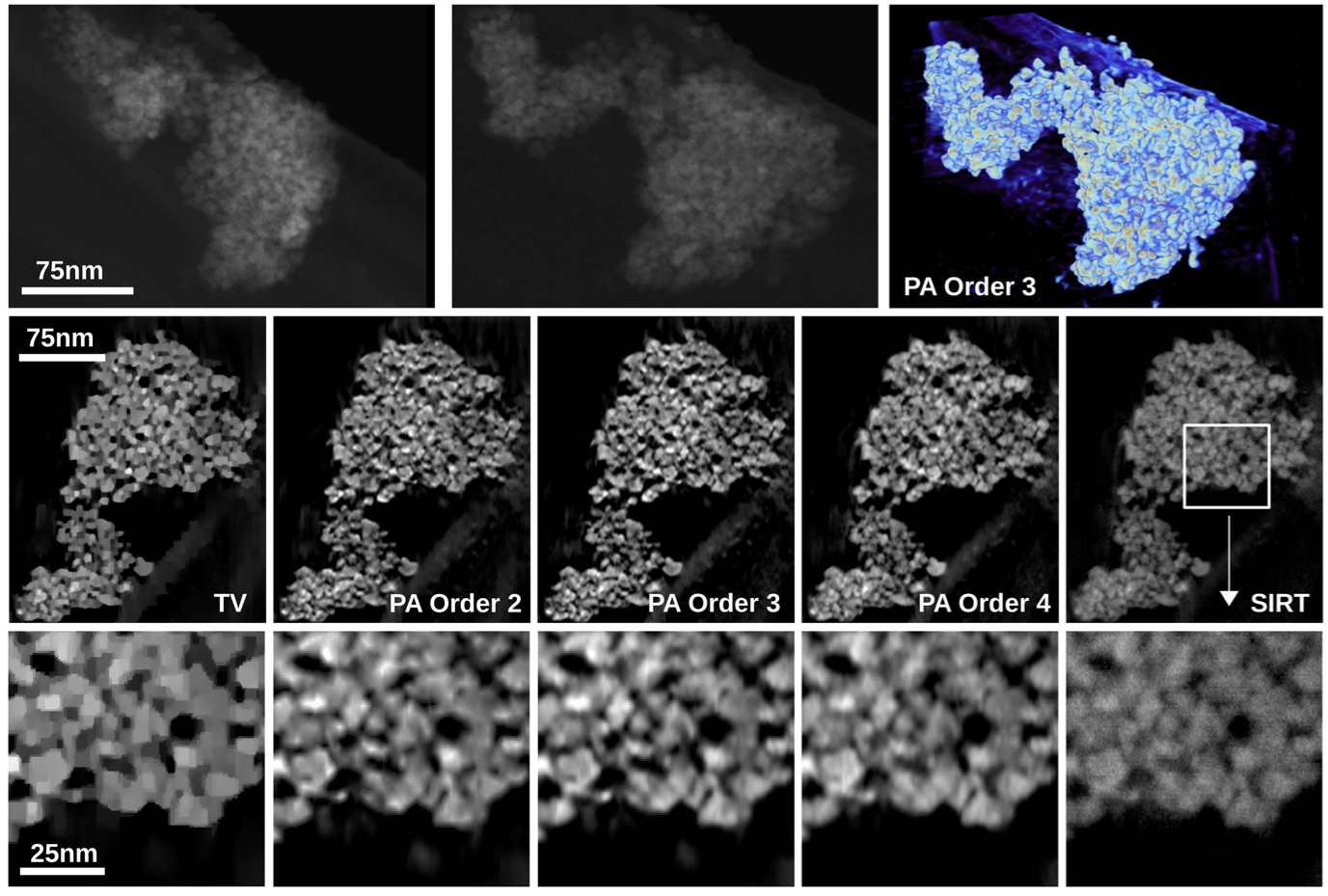}
 \caption{Top row: projection images of the nanoparticle taken at relative angles of (from left to right) $-46\degree$, $10\degree$, with the 3D volume rendering of the PA order 3 solution on the right.  Middle row: cross section image of the reconstructed volumes using (from left to right) TV, PA of order 2, 3, and 4, and least squares.  Bottom row: a magnified patch of each reconstruction.  The images demonstrate that PA is able to more accurately approximate the fine pore structure.}
 \label{bn-main}
\end{figure*}

\begin{figure*}
 \centering
 \includegraphics[scale = .65]{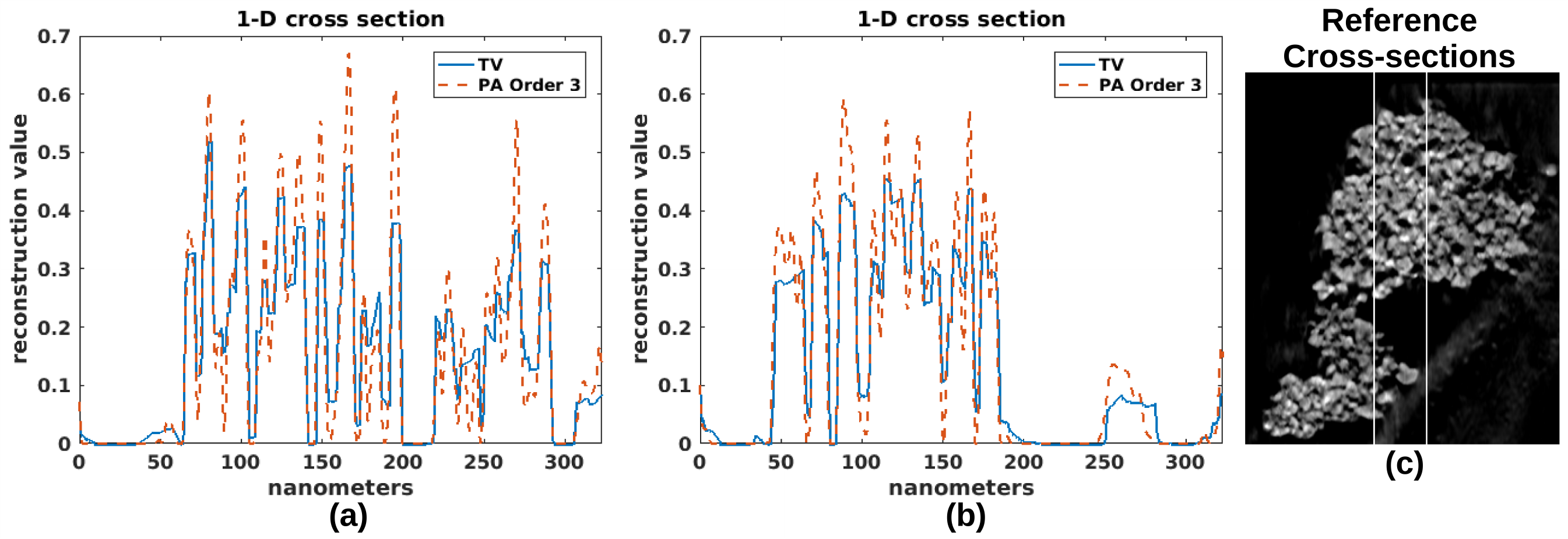}
 \caption{One-dimensional cross-sectional plots of the 3-D volume reconstructed with TV (blue, solid line) and PA order 3 (orange, dashed line), with the reference cross-sections shown in (c).  This plot allow us to see more clearly that the PA order 3 solution resolves the apparent pores much better than TV solution, where the solution from PA takes larger jumps both up and down.  This would ease the task of image segmentation.}
 \label{one-dimensional-plot}
\end{figure*}

\begin{figure*}
 \centering
 \includegraphics[scale=.34]{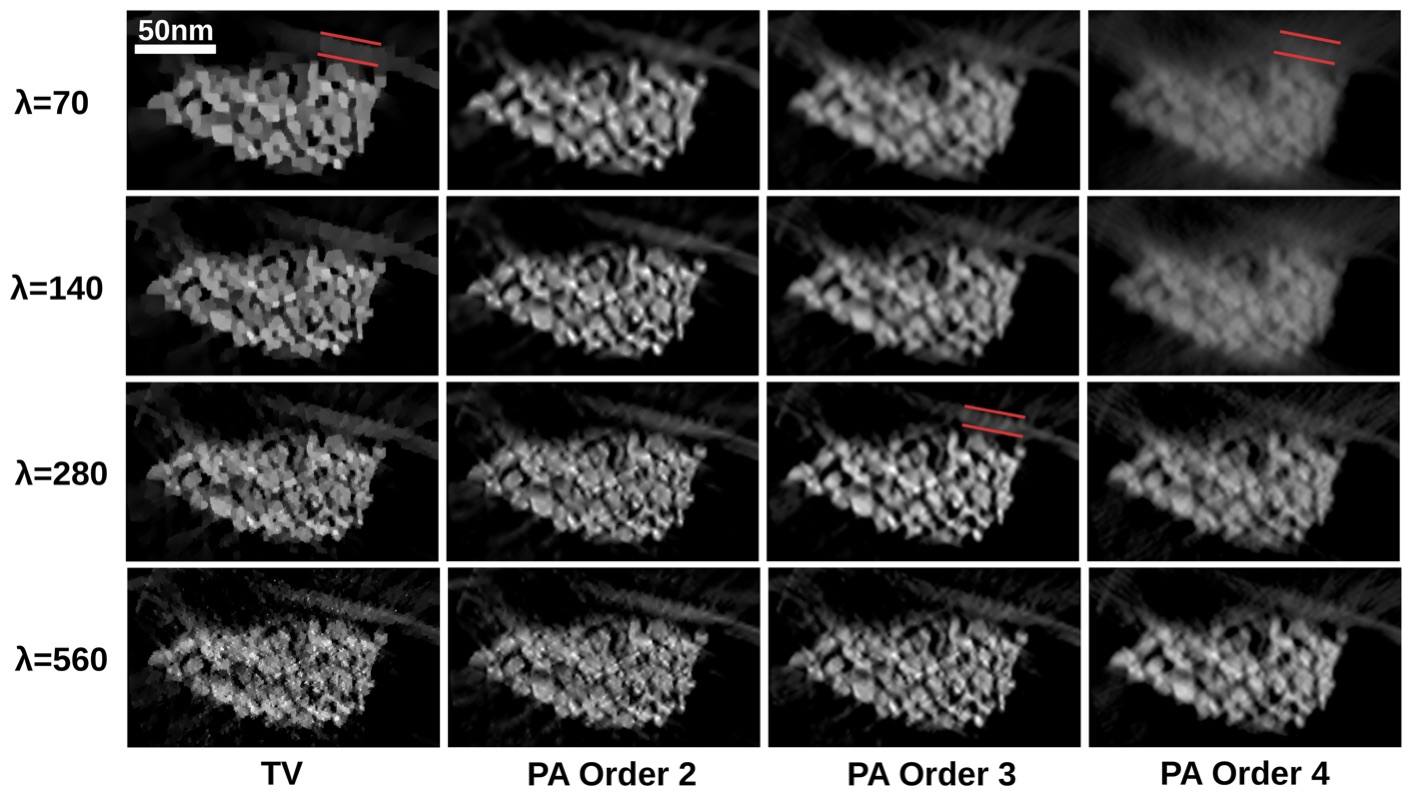}
 \caption{A reconstructed slice using (\ref{gen-PA}) for different regularization terms $T_k$ (columns) and different values of $\lambda$ (rows).  The suggested value for $\lambda$ for PA order $k$ is $\lambda = 2^{k-1}\cdot \lambda_1$, where $\lambda_1$ is the optimal choice for the TV formulation.  This means our suggested appropriate $\lambda$ selection for each column here is located on the main diagonal, which is verified by the quality of the images.  Red lines also indicate some missing wedge artifacts.}
 \label{bn-lambda}
\end{figure*}

The particular data set utilized in this investigation is of a porous, nanopolycrystalline, cubic boron nitride nanoparticle synthesized from periodic mesoporous hexagonal boron nitride at $10$ GPa and $1000\degree$C.  The electron microscopy was performed at $200$ kV on an FEI Tecnai instrument in STEM mode. The specimen grid was loaded into a Fischione single-tilt tomography holder.  The projections were acquired using the automated FEI acquisition software over a 144 angular range every $2$ degrees, giving a total of $73$ angles.  Each $1024\times 1024$ image was collected on a Fischione high angle annular dark field detector.  Two projections acquired at $-46\degree$ and $10\degree$ are shown (from left to right) in the top row of Figure \ref{bn-main}.  The ``missing wedge'' of $36\degree$ in this data set is the prototypical limit for ET imaging, and it is often times much larger.  We note that it is possible to obtain accurate results given less data \cite{DIPS}.  However, as our underlying image has many unstructured channels, we will utilize the full data set.

Prior to reconstruction, the projections were found to be accurately aligned using cross correlation, along with a manual search for the positioning of the axis of rotation. The full 3-D reconstructions were calculated using a least squares method (SIRT, \cite{SIRT}), TV regularization, and PA regularization of orders 2, 3, and 4.\footnote{{All associated MATLAB codes for the alignment and reconstruction algorithms are available at \cite{toby-web}.}}  After proper rescaling of the data values, we used an initial experiment to choose $\lambda_1=70$ for the TV formulation of (\ref{gen-PA}). We then applied (\ref{eq:lambdaselection}) for each increased value of $k$.  Finally, we also enforced a non-negative solution constraint $f\ge 0$ via a projected gradient method to decrease the size of the search space.  This is reasonable since $f$ is a nonnegative density function.

For image quality comparisons, the second row of Figure \ref{bn-main} shows a 2-D cross section slice of each reconstruction.  For closer inspection, the bottom row of the figure show a small magnified patch of the image.  This patch is indicated in second row of the rightmost image.  As a point of reference, observe that the standard least squares solution (SIRT) is noisy and does not capture the pores apparent in the other reconstructions.

The TV solution effectively eliminates the noise present in the least squares solution.  However, the lack of resolution causes the smaller pores to appear ``smeared out''.  In particular, while TV regularization encourages piecewise constant solutions, the data are too under-sampled to separate the small pores.  We note that given sufficient sampling and/ or spacing of the pores, TV regularization would be ideal for this kind of piecewise constant solution.

Put another way, an accurate depiction of small pores, each consisting neighboring jumps in opposing directions (see Figure \ref{one-dimensional-plot}), increases the value of the $\ell_1$ penalty term.  However the contribution to the corresponding fidelity term, which accounts for data fit in the approximation, is not significant since the pores are small and do not constitute a large part of the overall data fitting.  In general, as the size of the pores increase, how accurate their depiction is becomes more significant in the fidelity term, as the contribution constitutes a larger part of the overall approximation.  On the other hand, if a pore is not detected, then the value of the $\ell_1$ penalty term is small, but the corresponding fidelity term will increase due to the poor data fit.  Thus if the problem is well resolved, using TV regularization provides the best solution, as the piecewise constant nature of the solution will provide a good data fit for the fidelity term and the sparsity of the pore edges is captured in the $\ell_1$ term.

In this investigation we are concerned with {\em under-resolved} data, and in this case, the TV regularization essentially ``smears'' neighboring pixel values around small pores to reduce the size of the penalty term.  However, it is evident from the least squares solution (SIRT) in Figure \ref{bn-main} (far right) that some regularization is needed to remove unwanted artifacts.  Using the PA penalty term of orders $2$ or $3$ appear to accurately capture the fine porous details, while suppressing the noise.  Specifically, the pores are captured with edges which are closer to first and second degree polynomials.  Using PA of order $4$, which encourages piecewise polynomials of degree less than four, allows for an approximation in which the noise once again becomes apparent in the image.  

The 1-D cross-sectional plots of the reconstructions in Figure \ref{one-dimensional-plot} illustrate with additional clarity why using PA (shown for $k =3$) is more effective than TV.  While both plots show similar jump locations, TV is unable to capture full variation apparent in the PA solution.  The TV solution effectively smears over the image values where the PA solution captures much larger jump values via polynomial fit.

{The penalty term influence can be reduced by increasing $\lambda$.  However, this will not improve the approximation of the small pores since noise in the data will cause a poor fidelity fit.  Figure \ref{bn-lambda} illustrates this issue by showing a slice of the reconstructed volume using several choices of regularization parameter $\lambda$.  Observe also that along each diagonal band of this figure, $\lambda$ is chosen using (\ref{eq:lambdaselection}). The clear degradation of the solution off these bands demonstrates the validity of our procedure for choosing $\lambda$ robustly.}

Figure \ref{bn-lambda} also gives us a view of the reconstruction in the direction of the missing wedge.  In this case, we would observe vertical elongation of these images due to the missing wedge.  A careful examination of the elements along the diagonal band do show reduced elongation for orders 2 and 3 compared with TV, which again can be characterized by ``smearing'' of the pores.  Also observe for example, the carbon grid on which the particle rests (positioned in the northeast corner of the images) exhibits reduced elongation with the higher order methods.  This is indicated with the red lines outlining this carbon support in several of the images. However, the same artifacts are observable in the horizontal direction as well to a somewhat lesser extent. 

\FloatBarrier

\section{Comparison with Discrete Tomography}

\begin{figure*}
 \centering
 \includegraphics[width=.9\textwidth]{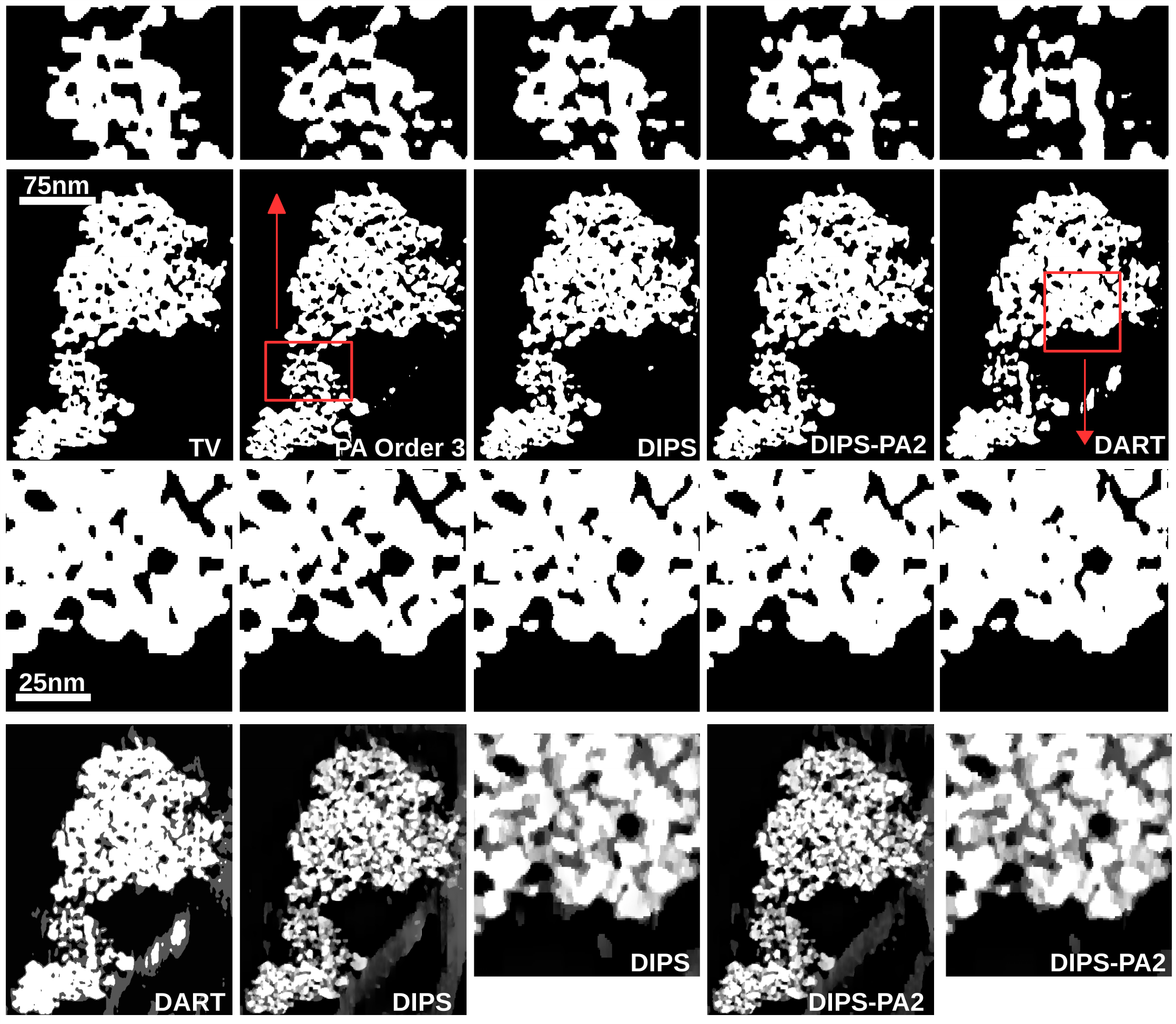}
 \caption{Segmented TV and PA solutions compared with discrete tomography from DIPS, DIPS-PA, and DART.  The top and third rows are magnified patches of the second row, with their location indicated by the red boxes and arrow.  The bottom row shows the discrete tomographic reconstructions before segmenting.  The DIPS solutions and segmented PA order 3 appear the most consistent.  The segmented TV misses some of the pore structure, and DART shows some particle morphology that disagrees with all other solutions.}
 \label{bn-DT}
\end{figure*}

Scientific quantification of reconstructed nanoparticles in electron tomography usually requires an accurate segmentation of the reconstructions, which can sometimes reduce to a manual task biased by the human error. To reduce human error and missing wedge artifacts, methods have recently emerged that implement qualitative segmentation techniques \emph{within} the reconstruction algorithm \cite{DART,DIPS}.  These methods are known as discrete tomography \cite{hermanDT}, lending its name to the discrete gray levels corresponding to the various a priori assumed densities of the materials within a given sample. Given several gray levels $\mathbf{\rho} = \{\rho_1,\rho_2,\dots, \rho_k\}$, where the number of gray levels, $k$, is typically 2 or 3, the general discrete tomography problem is to recover
\begin{equation}\label{DTform}
 f = \mathop{\mathrm{argmin}}_{f:~ f_i\in \rho} \| Wf - b\|_2^2.
\end{equation}

While the prior knowledge of exact gray levels in the reconstruction may seem advantageous, we mention a few inherent drawbacks:
\begin{itemize}
\item Solving the general formulation of discrete tomography is NP hard, hence a special algorithm must be carefully designed in hopes of only approximating the solution to the general formulation. 
 \item In practice, the gray levels must be accurately determined to achieve good results.  Streamline determination of the gray levels can be a challenging barrier to many working scientists, and one that is not presented with more general $\ell_1$ algorithms.
 \item The tight constraint that the reconstruction must contain only several gray values leaves it highly susceptible to noisy data.  Moreover, if the assumption of a few gray levels is incorrectly prescribed, then the resulting reconstruction may be highly inaccurate. 
\end{itemize}

Nevertheless, such algorithms have been shown to yield accurate results in numerous instances \cite{GorisTVDART,Jinschek2008589,Batenburg2009730}.  For comparison, we present results on popular discrete algebraic reconstruction technique (DART) developed and extensively studied by Batenburg and others \cite{DART}, as well as the discrete iterative partial segmentation technique (DIPS) developed in \cite{DIPS}.  Since this is not the main focus of this paper, for a detailed description of these methods the reader should look into the referenced papers.

Both techniques work through a sequence of segmentation and refinement steps.  DART iterates between segmenting the solution, and refining after each segmentation with a smoothing procedure (usually through convolution with a kernel) and a refinement along the boundaries of the image based on the error between the segmented solution and the data.  DIPS takes a more refined approach by combining $\ell_1$ optimization with discrete tomography.  We note similar ideas are developed in \cite{TVR-DART}.  Given a current solution, a partial segmentation is implemented in which only the pixels that have high probability of falling into the gray levels set are segmented, based of the the distances from the gray levels.  This solution is then refined over the subset of unsegmented pixels, a problem which has a significantly reduced dimension, and this refinement is carried out by implementing an augmented TV regularized optimization functional.  The steps of partial segmentation and refinement with the TV functional are iterated as many times as necessary or chosen.  Although the final solution may not be completely discrete, the idea is that the solution accuracy increases with each iterate as the problem dimension is reduced with a greater number of pixels becoming classified.

Naturally, an additional contribution of this work is to modify the DIPS approach to use an augmented PA regularized optimization functional in place of the original TV functional, which we have included in our results.  In particular, in equation (17) of \cite{DIPS}, $\|f\|_{TV}$ is replace with $\|T_kf \|_1$, where here we simply use $k=2$.

Comparisons with discrete tomography on the BN nanoparticle are shown in Figure \ref{bn-DT}, where we again show the same slice from Figure \ref{bn-main} for thorough comparison.  In the top and bottoms rows of the figure are magnified smaller patches indicated by the red boxes and arrows.  The TV, PA, and SIRT solutions were segmented with a fixed threshold value at 0.15, one of the thresholds used in the DART reconstruction.  The DIPS and DART reconstructions used three gray levels, one for the background, the carbon grid (removed from the images), and the BN particle.  

Reconstructions from DIPS and PA show the greatest similarity, most compatible with what was observed in the continuous solutions.  As expected from the continuous results, the segmented TV solution misses some of the fine details.  The DART solution maintained some level of accuracy, as evidenced by the image in the bottom row.  Although the magnified patch for DART in the top row shows significant discrepancy from all of the others.  Very little differences exists between the original DIPS and the modified DIPS using PA (labeled DIPS-PA2).  It seems though DIPS-PA2 has slightly more pore structure than the original DIPS, which would agree with the general findings of this paper.

\section{Simulations Results}
\begin{figure*}
 \centering
 \includegraphics[scale=.35]{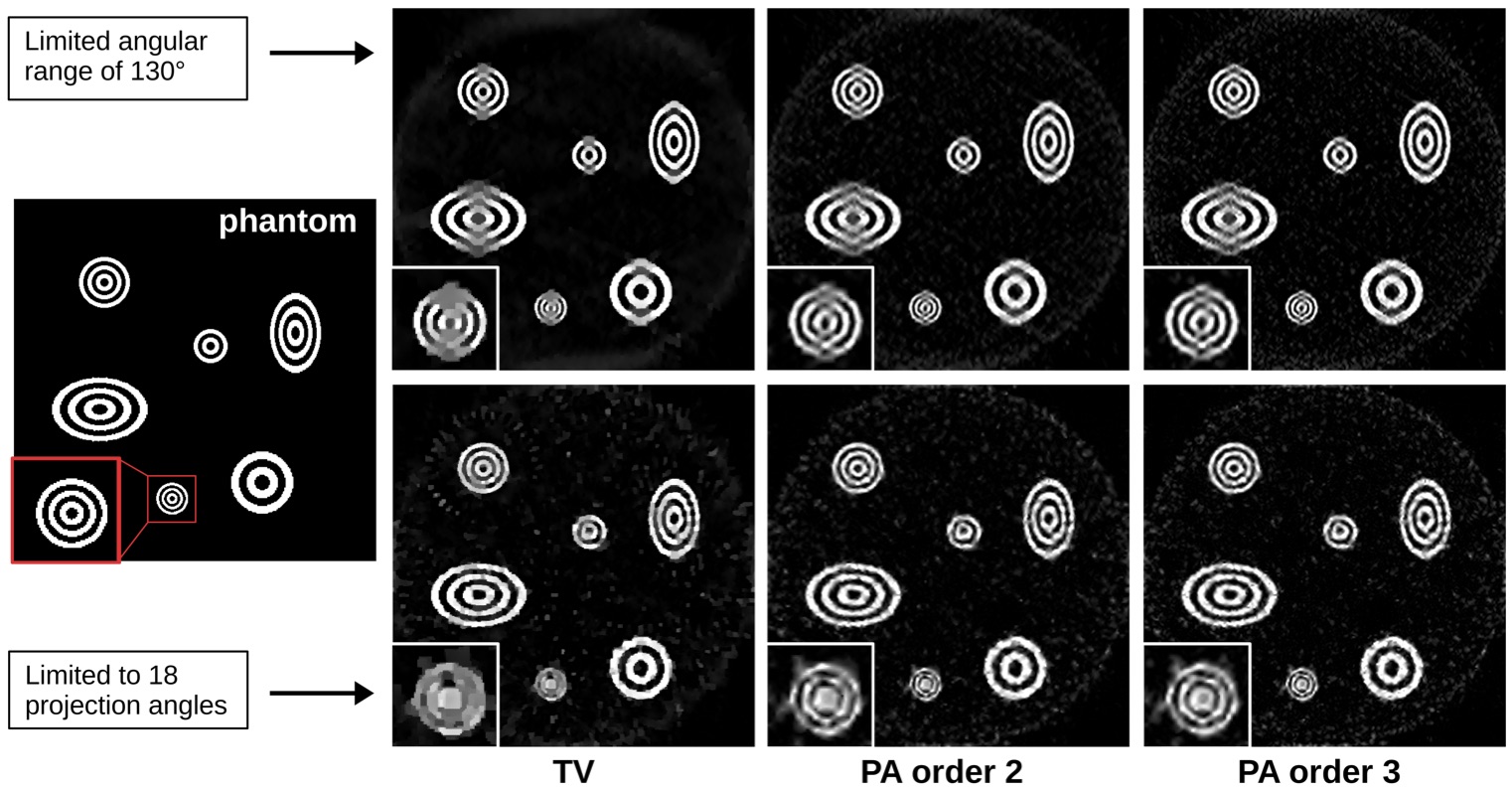}
 \caption{Results from tomography simulations.  The smallest set of concentric rings is magnified in the bottom left of each image.}
 \label{circle-fig}
\end{figure*}

\begin{figure*}
 \centering
 \includegraphics[scale=.7]{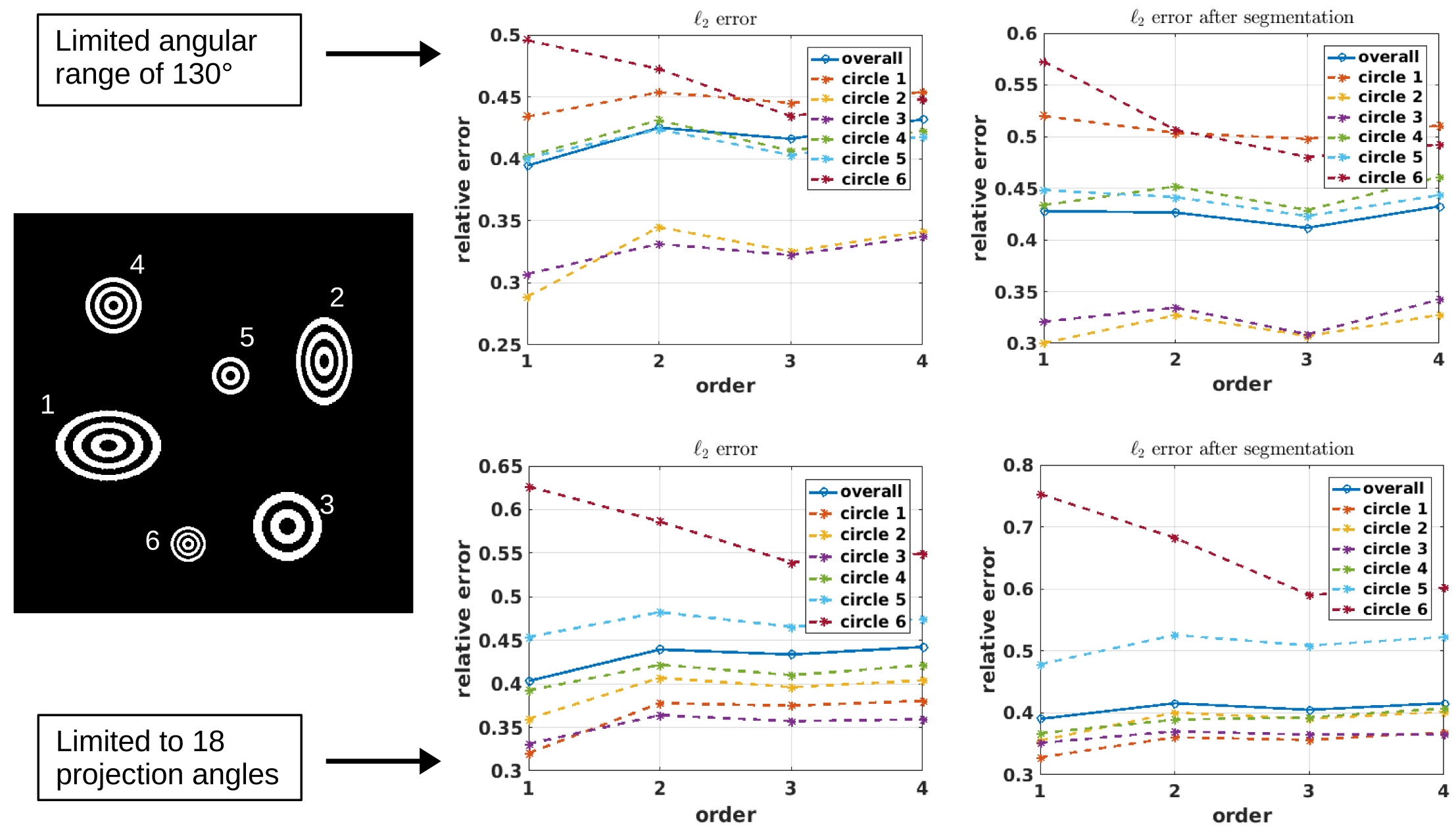}
 \caption{The relative $\ell_2$ error measures from tomographic simulations as a function of the order of the regularization used for the reconstruction.  The global error is given as well as the local error for each set of concentric circles, which are numbered as indicated in the left panel.  The error is also measured after segmentation is applied (right column).}
 \label{circle-fig-error}
\end{figure*}
To further analyze (\ref{gen-PA}) for different regularization terms, we consider a test phantom image, shown in Figure \ref{circle-fig}.  The image is binary and contains a series of concentric rings of various sizes and spacings.   As discussed previously, TV regularization is well suited for piecewise constant images when the data are sufficiently sampled.  However, the details in this image given under-resolved data make the reconstruction more challenging.

We generate the data with Poisson noise with fixed mean  to simulate noise similar to that in electron microscopy \cite{Mevenkamp}, and then consider the following two cases: 
\begin{itemize}
\item {\bf Missing wedge:} Here we assume the projection data are collected at only $130\degree$ out of the possible $180\degree$.  The projections are generated at $\Delta\theta = 2.5\degree$ within the $130\degree$ angular range.  
\item {\bf Limited data:} Here we assume we have the full $180\degree$ range.  The angle step, $\Delta\theta = 10\degree$, leaving us with only $18$ projections. 
\end{itemize}

As mentioned previously, the added noise and under sampling suggest the use of TV regularization, which is very effective for recovering piecewise constant functions from under sampled data, assuming certain properties of the sampling distribution.  For example, random sampling throughout the domain would suffice.  However, as described above, the data sampled in our tomography  simulation are more restricted.  Nonetheless,  Figure \ref{circle-fig} demonstrates that employing convex optimization in the form (\ref{gen-PA}) using the PA regularization of orders $2$ or $3$ is effective in reconstructing the image. For clarity, the smallest of the concentric rings is magnified in the bottom of left of each image.  Observe that for the missing wedge case (top row), the missing angles cause some of the concentric rings to appear ``smeared'' together in the vertical direction, an effect that is more pronounced in the TV solution.  These findings agree with the results for the experimental ET data sets, where some of the fine features also appeared to be smeared together. While the PA regularization causes less smearing within the individual features, some artifacts become apparent away from the rings.  One explanation for the differences in these reconstructions may be that while the TV regularized solution distributes the noise throughout the reconstruction, the PA regularized solution causes noise to appear more ``point-like'', an effect visible in Figure \ref{one-dimensional-plot}.  Therefore, post-processing with a threshold to reduce the noise is more straightforward with TV.  Post-processing algorithms to reduce the effects of noise when using PA regularization are discussed in \cite{phd-denker}.

For the reduced resolution case, shown in the bottom row of Figure \ref{circle-fig}, a similar effect is seen in that some of the smaller concentric rings are again smeared together in the TV regularized solution, while PA regularization yields more details.  However, since the larger rings are sufficiently resolved, both the TV and PA regularization algorithms accurately capture their details.  The noise artifacts away from the rings are point-like in all cases, but again appear to be more evenly distributed in the TV regularized solution.

The relative $\ell_2$ error measures are provided in Figure \ref{circle-fig-error}.  We provide the global and local error measures for each of the 6 sets of concentric circles, The numbering of these sets of circles is indicated in the image to the left, and the ordering is approximately from largest to smallest.  For the larger sets of circles, the TV solution (order 1) and PA order 3 solutions show modest improvement over order 2, and TV appears marginally better than order 3.  Generally speaking, for the larger sets of circles each reconstruction recovered the structure well, and the small differences in the error is primarily due the noise.  For the set number 6 however, the higher orders show notable improvements over TV, and these gains begin to diminish after order 3.  The plots in the rightmost column show the errors after segmentation of the solutions, where for simplicity a simple threshold was used to transform the images into their binary segmentations.  The plots after segmentation generally follow the same trends as before segmentation, perhaps indicating slightly more improvements for the higher orders.

\section{Concluding Remarks} 
\label{sec:conclusion}
This investigation demonstrates that applying higher order TV regularization, namely the polynomial annihilation, to under-resolved electron tomography data yields accurate results in the sense that small pores and channels in the images can be captured and the noise can be adequately suppressed.  In particular, the PA regularization appears to be more effective at capturing the fine details of the image both for under-resolved data and the missing wedge case.  Our comparisons with discrete tomography suggest that DIPS reconstructions also yield accurate results comparable with higher order PA.  The challenges associated with experimental discrete tomography, namely gray value determination and susceptibility to noise as pointed out in Figure \ref{bn-DT} for DART, still makes it somewhat limited.  The findings from experimental results with PA were additionally verified in tomographic simulations. We also provided numerical evidence for our selection of the regularization parameter for each increasing order of PA regularization.  

\bibliographystyle{abbrv}

\end{document}